\documentclass{amsproc}
\usepackage{eurosym}
\usepackage{amssymb}
\usepackage{amsfonts}
\usepackage{cmtt}

\theoremstyle{plain}
\newtheorem{theorem}{Theorem}
\newtheorem{corollary}[theorem]{Corollary}

\newtheorem{question}[theorem]{Question}
\numberwithin{equation}{section}
\DeclareMathOperator{\ch}{ch}

\begin{document}

\title{The coloring game on matroids}

\author[Micha\l\ Laso\'{n}]{Micha\l\ Laso\'{n}}

\dedicatory{\upshape
Institute of Mathematics of the Polish Academy of Sciences,\\ ul.\'{S}niadeckich 8, 00-656 Warszawa, Poland\\ \textmtt{michalason@gmail.com}
\smallskip
\\ Theoretical Computer Science Department,\\ Faculty of Mathematics and
Computer Science, Jagiellonian University,\\ ul.{\L}ojasiewicza 6, 30-348 Krak\'{o}w, Poland
\\ \textmtt{mlason@tcs.uj.edu.pl}}

\thanks{Research partially supported by the Polish National Science Centre grant no. 2011/03/N/ST1/02918. The paper was completed during author's stay at Freie Universit\"at Berlin in the frame of Polish Ministry ``Mobilno\'s\'c Plus'' program.}
\keywords{Matroid coloring, Acyclic coloring, Game coloring, Game arboricity.}

\begin{abstract}
A coloring of the ground set of a matroid is \emph{proper} if elements of the same color form an independent set. For a loopless matroid $M$, its \emph{chromatic number} $\chi (M)$ is the minimum number of colors in a proper coloring. In this note we study a game-theoretic variant of this parameter. 

Suppose that Alice and Bob alternately properly color the ground set of a matroid $M$ using a fixed set of colors. The game ends when the whole matroid has been colored, or if they arrive to a partial coloring that cannot be further properly extended. Alice wins in the first case, while Bob in the second. The \emph{game chromatic number} of $M$, denoted by $\chi_{g}(M)$, is the minimum size of the set of colors for which Alice has a winning strategy. Clearly, $\chi_{g}(M)\geq\chi (M)$. 

We prove an upper bound $\chi_{g}(M)\leq 2\chi (M)$ for every matroid $M$. This improves and extends a result of Bartnicki, Grytczuk and Kierstead \cite{BaGrKi08}, who showed that $\chi_{g}(M)\leq 3\chi (M)$ holds for graphic matroids. Our bound is almost tight, as we construct a family of matroids $M_k$ (for $k\geq 3$) satisfying $\chi (M_k)=k$ and $\chi_{g}(M_k)=2k-1$.
\end{abstract}

\maketitle

\section{Introduction}

Let $M$ be a matroid on a ground set $E$ (we refer the reader to \cite{Ox92} for a background of matroid theory). A \emph{coloring} of $M$ is an assignment of colors to the elements of $E$. In analogy to graph theory we say that a coloring is \emph{proper} if elements of the same color form an independent set in the matroid. Via this correspondence one can define for matroids all chromatic parameters studied for graphs. 

The \emph{chromatic number} of a loopless matroid $M$, denoted by $\chi(M)$, is the minimum number of colors in a proper coloring of $M$. For instance, if $M$ is a graphic matroid obtained from a graph $G$, then $\chi(M)$ is the least number of colors needed to color edges of $G$ so that no cycle is monochromatic. This number is known as the \emph{arboricity} of the underlaying graph $G$. 

The chromatic number of a matroid can be easily expressed in terms of its rank function. Edmonds \cite{Ed65} gave an explicit formula, extending a theorem of Nash-Williams \cite{Na64} for graph arboricity. Seymour \cite{Se98} (see also \cite{La13}) proved that the chromatic number $\chi(M)$ of a matroid $M$ is equal to the list chromatic number $\ch(M)$ -- another well-known concept from graph theory, initiated by Vizing \cite{Vi76}, and independently by Erd\H{o}s, Rubin and Taylor \cite{ErRuTa80}.

In this note we study a game-theoretic variant of the chromatic number of a matroid (for other variants see ex. \cite{La14,LaLu12}). It is defined by a game beetwen two players -- Alice and Bob. They alternately color elements of the ground set $E$ of a matroid $M$ using a fixed set of colors $C$. The only rule that both players have to obey is that at any moment of the play the partial coloring is proper. Who starts does not influence our results, but to make the definition strict suppose that Alice makes the first move. The game ends when the whole matroid has been colored, or if they arrive to a partial coloring that cannot be further extended (what happens when trying to color any uncolored element with any possible color results in a monochromatic circuit of $M$). Alice wins in the first case, while Bob in the second. The \emph{game chromatic number} of a matroid $M$, denoted by $\chi_{g}(M)$, is the minimum size of the set of colors $C$ for which Alice has a winning strategy. 

The above game is a matroidal analog of a well-studied graph game coloring, which was introduced by Brams (cf. \cite{Ga81}) for planar graphs with a motivation to give an easier proof of the four color theorem. The game was independently reinvented by Bodlaender \cite{Bo91}, and since then the topic developed into several directions leading to interesting results and challenging open problems (see a recent survey \cite{BaGrKiZh07}).

The first step in studying game chromatic number of matroids was made by Bartnicki, Grytczuk, and Kierstead \cite{BaGrKi08}. They showed that for every graphic matroid $M$ inequality $\chi_{g}(M)\leq 3\chi(M)$ holds. In Theorem \ref{main1} we improve and extend this result by proving that for every loopless matroid $M$ we have $\chi_{g}(M)\leq 2\chi (M)$. This gives a nearly tight bound, since in Theorem \ref{main2} we provide a family of matroids $M_k$ (for $k\geq 3$) satisfying $\chi(M_k)=k$ and $\chi_{g}(M_k)=2k-1$. Our bounds remain true also for the fractional parameters, as well as for list version of the game chromatic number. 

\section{A strategy for Alice}

To achieve our goal we shall need a more general version of the matroid coloring game. Let $M_{1},\dots ,M_{d}$ be matroids on the same ground set $E$, and let $\{1,\dots,d\}$ be the set of colors. As before, the players alternately color elements of $E$, but now for all $i$ the set of elements colored with $i$ must be independent in the matroid $M_{i}$. As before, Alice wins if at the end of the game the whole set $E$ is colored, Bob wins otherwise. We call this game \emph{coloring game on} $M_1,\dots,M_d$. The initial game on $M$ with $d$ colors coincides with the coloring game on $M_{1}=\dots=M_{d}=M$.

\begin{theorem}\label{main1}
Let $M_{1},\dots ,M_{d}$ be matroids on a ground set $E$. If there exist sets $V_{1},\dots,V_{d}$, with $V_{i}$ independent in $M_{i}$, such that $V_{1}\cup\dots\cup V_{d}=E\cup E$ as multisets, then Alice has a winning strategy in the coloring game on $M_{1},\dots ,M_{d}$. In particular, for every loopless matroid $M$ we have $\chi_{g}(M)\leq 2\chi (M)$. 
\end{theorem}

\begin{proof}
Let us fix a $2$-covering of $E$ by sets $V_1,\dots,V_d$ independent in corresponding matroids. Let $U_{i}$ be the set of elements colored (at a fixed moment of the play) with $i$. Then $C=U_{1}\cup\dots\cup U_{d}$ is the set of all colored elements. Alice will try to keep the following invariant after each of her moves:
\begin{equation}\label{cond}
\text{for each }i\text{, the set }U_{i}\cup (V_{i}\setminus C)\text{ is independent in }M_{i}.
\end{equation}
Moreover, element $e\in V_{i}\cap V_j$ will be colored by Alice only with $i$ or $j$.

Observe first that if the condition $(\ref{cond})$ holds and there is an uncolored element $e$ ($e\in V_i$ for some $i$), then the player can always make an `obvious' move, namely to color $e$ with $i$. After this move the condition $(\ref{cond})$ remains true.

To prove that Alice can keep the condition $(\ref{cond})$ assume that the condition held after her previous move (let $U_i$ and $C$ be defined at this moment) and then Bob colored an element $e$ with color $j$. 

If $(U_{j}\cup e)\cup (V_{j}\setminus C)$ is independent in $M_{j}$, then $(\ref{cond})$ still holds and we use the above observation. 

When $(U_{j}\cup e)\cup (V_{j}\setminus C)$ is dependent in $M_{j}$, then by the 
augmentation axiom we can extend the independent set $U_{j}\cup e$ from the independent set $U_{j}\cup (V_{j}\setminus C)$ in $M_{j}$. The extension equals to $(U_{j}\cup e)\cup (V_{j}\setminus C)\setminus f$ for some $f\in V_{j}\setminus C$. Since sets $V_{1},\dots,V_{d}$ form a $2$-covering we know that $f\in V_{l}$ for some $l\neq j$. Now Alice has an admissible move, and her strategy is to color $f$ with color $l$. It is easy to observe that after her move the condition $(\ref{cond})$ is preserved.

To get the second part of the assertion suppose that $\chi(M)=k$. By the first part of the assertion applied to $M_{1}=\dots=M_{2k}=M$ we infer that Alice has a winning strategy with $2k$ colors, and therefore $\chi_{g}(M)\leq 2k$. 
\end{proof}

As usually in this kind of games the following question seems to be natural and non trivial (for graph coloring game it was asked by Zhu \cite{Zh99}).

\begin{question}
Suppose Alice has a winning strategy in the coloring game on a matroid $M$ with $k$ colors. Does she also have a winning strategy with $l>k$ colors?
\end{question}

We generalize our result to list version. The rules of the game between Alice and Bob do not change, except that now each element has its own list of available colors. So, instead of a fixed set of colors $C$, each element $e\in E$ can be colored by Alice or Bob only with a color from its list $L(e)$. The minimum number $k$ for which Alice has a winning strategy for every assignment of lists of size $k$ is called the \emph{game list chromatic number} of $M$, and denoted by $\ch_{g}(M)$. Clearly, $\ch_{g}(M)\geq\chi_{g}(M)$ for every loopless matroid $M$, however the same upper bound as for $\chi_{g}(M)$ works also for list parameter.

\begin{corollary}
For every loopless matroid $M$ we have $\ch_{g}(M)\leq 2\chi(M)$.
\end{corollary}

\begin{proof}
Let $L$ be a list assignment of size $2\chi(M)$. Without loss of generality we can assume that lists $L(e)$ are contained in the set $\{1,\dots,d\}$ for some $d$. By a theorem of Seymour \cite{Se98}, there is a $2$-covering of $E$ by sets $V_1,\dots,V_d$, such that elements of $V_i$ have color $i$ on their lists and $V_i$ is independent in $M$. Let $M_i$ be the matroid $M$ restricted to the set $V_i$ with the ground set trivially extended to $E$. By Theorem \ref{main1} Alice has a winning strategy in the coloring game on lists $L$.
\end{proof}

In view of Seymour's Theorem \cite{Se98} we ask the following natural question.

\begin{question}
Does for every matroid $M$ the game list chromatic number equals to the game
chromatic number $\ch_g(M)=\chi_g(M)$?
\end{question}

\section{A strategy for Bob}

We present a family of transversal matroids $M_{k}$ for $k\geq 3$, with $\chi(M_k)=k$ and $\chi_{g}(M_k)\geq 2k-1$. This slightly improves the lower bound from the paper of Bartnicki, Grytczuk and Kierstead \cite{BaGrKi08}. They gave an example of graphic matroids $H_{k}$ satisfying $\chi(H_{k})=k$ and $\chi_{g}(H_{k})\geq 2k-2$ for every $k\geq 1$. 

Fix $k\geq 3$. Let $D_1,\dots,D_{3k(2k-1)}$ be disjoint sets such that each set $D_i=\{d_{1,i},\dots,d_{k,i}\}$ has exactly $k$ elements. Let $C=\{c_{1,1},\dots,c_{k,2k-1}\}$ be a set with $k(2k-1)$ elements, disjoint from sets $D_i$. Let $E$ be the union of sets $C$ and $D_1,\dots,D_{3k(2k-1)}$. Let $M_k$ be the transversal matroid (see \cite{Ox92} for a definition) on a ground set $E$ with multiset of subsets $\mathcal{A}$ consisting of sets $D_1,\dots,D_{3k(2k-1)}$ and $(2k-1)$ copies of $E$. In other words, a subset $I\subset E$ is independent in $M_k$ if there is some subset $J\subset I$ with $\left\vert J\right\vert\leq 2k-1$, such that $I\setminus J$ contains no elements of $C$ and at most one element of each $D_i$.

\begin{theorem}\label{main2}
For $k\geq 3$ matroid $M_{k}$ satisfies $\chi(M_k)=k$ and $\chi_{g}(M_{k})\geq 2k-1$.
\end{theorem}

\begin{proof}
To prove the first part observe that we can partition the ground set $E$ into $k$ independent sets $V_i=\{c_{i,1},\dots,c_{i,2k-1},d_{i,1},\dots,d_{i,3k(2k-1)}\}$. 

To prove the second part notice that the rank of $C\cup D_i$ equals to $2k$, so if there are $t$ elements from $D_i$ colored with $i$, then there are at most $2k-t$ elements from $C$ colored with $i$. This suggests that Bob should try to color each $D_i$ with one color. 

Suppose $\{1,\dots,h\}$ is the set of colors in the game, and $h\leq 2k-2$. We will describe a winning strategy for Bob. Alice will always loose the game because she will be not able to color all elements of the set $C$. 

Assume first that Alice colors only elements from the set $C$, and her goal is to color all of them. It is the main case to understand. Denote by $d_i$ and $c_i$ the number of elements colored with $i$ in $D_i$ and $C$ respectively (at some moment of the play). Bob wants to keep the following invariant after each of his moves:
\begin{equation}\label{bob1}
d_i\geq c_i\text{, for every color }i.
\end{equation}
It is easy to see that he can always do it, and in fact there is an equality $d_i=c_i$ for every $i$. Bob just mimics Alice's moves. Whenever she colors some $c\in C$ with $i$ he responds by coloring an element of $D_i$ with $i$. He can do it, because when after Alice's move $c_i=d_i+1 $ for some $i$, then $c_i+d_i\leq 2k$. But, then also $c_i+(d_i+1)\leq 2k$, so elements colored with $i$ are independent, and $d_i+1\leq k$, so there was an uncolored element in $D_i$. 

Observe that when Bob plays with this strategy, then for each color $i$ we have $c_i+d_i\leq 2k$ and $c_i\leq d_i$, so as a consequence $c_i\leq k$. This means that Alice can color only $hk$ elements of $C$, thus she looses.

It remains to justify that coloring elements of $D_1\cup\dots\cup D_{3k(2k-1)}$ by Alice can not help her in coloring elements of $C$. To see this we have to modify the invariant that Bob wants to keep. We assume $k\geq 4$, because for $k=3$ more careful case analysis is needed. Denote by $d_i,f_i$ the number of elements in $D_i\cup D_{i+h}\cup D_{i+2h}$ colored with $i$ and with some other color respectively, and $c_i$ as before. Now the invariant that Bob wants to keep after each of his moves is the following: 
\begin{equation}\label{bob2}
d_i\geq c_i+f_i\text{ or }d_i\geq k+2\text{, for every color }i.
\end{equation}
Analogously to the previous case one can show that Bob can keep this invariant. He just have to obey one more rule. Whenever Alice colors an element of $D_i\cup D_{i+h}\cup D_{i+2h}$, then Bob colors another element of this set with $i$ always trying to keep $\epsilon_i$, the number of sets among $D_i,D_{i+h},D_{i+2h}$ with at least one element colored with $i$, as low as possible. This completes the description of Bob's strategy. Notice that $f_i\geq\epsilon_i-1$.

Condition $(\ref{bob2})$ gives the same consequence as $(\ref{bob1})$. Indeed, let $\mathcal{D}_i$ be a the union of those of $D_i,D_{i+h},D_{i+2h}$, which have an element colored with $i$. Matroid $M_k$ restricted  to the set $C\cup\mathcal{D}_i$ has rank $2k-1+\epsilon_i$. If $d_i\geq c_i+f_i$, then 
$$2k-1+\epsilon_i\geq  c_i+d_i\geq 2c_i+f_i\geq 2c_i+\epsilon_i-1.$$ 
Otherwise, if $d_i\geq k+2$, then
$$2k+2\geq 2k-1+\epsilon_i\geq  c_i+d_i\geq c_i+k+2.$$ 
In both cases $c_i\leq k$, so there can be at most $k$ elements in $C$ colored with $i$. 
\end{proof}

Our results lead naturally to the following question.

\begin{question}
Does the inequality $\chi_g(M)\leq 2\chi(M)-1$ hold for every loopless matroid $M$?
\end{question}

\section*{Acknowledgements}

I would like to thank Jarek Grytczuk for many inspiring conversations, in particular for introducing me to game arboricity and posing a question for arbitrary matroids. Additionally, I want to thank him for the help in
preparation of this manuscript.


\end{document}